\documentclass[letterpaper, 10 pt, conference]{ieeeconf} 
\usepackage[latin9]{inputenc}
\usepackage{float}
\usepackage{amsmath}
\usepackage{amssymb}
\usepackage{graphicx}
\usepackage{todonotes}
\usepackage{gensymb}
\usepackage{booktabs} 
\usepackage{algorithm,algpseudocode} 

\IEEEoverridecommandlockouts                              
\global\long\def\norm#1{\lVert#1\rVert}%
\title{\textbf{Guided Policy Search using Sequential Convex Programming for Initialization of  Trajectory Optimization Algorithms}}

\author{Taewan Kim, Purnanand Elango, Danylo Malyuta, and 
Beh\c{c}et A\c{c}\i kme\c{s}e\thanks{Autonomous
Controls Laboratory, William E. Boeing Department of Aeronautics and Astronautics, University
of Washington, Seattle, WA 98105, USA \{twankim, pelango, danylo,
behcet\}@uw.edu. This work is supported in part by Office of Navel Research under Grant N00014-20-1-2288 and  AFOSR Grant FA9550-20-1-0053. }}
\begin{document}
\maketitle
\thispagestyle{empty}
\pagestyle{empty}
\begin{abstract}
Nonlinear trajectory optimization algorithms have been developed to handle optimal control problems with nonlinear dynamics and nonconvex constraints in trajectory planning. The performance and computational efficiency of many trajectory optimization methods are sensitive to the initial guess, i.e., the trajectory guess needed by the recursive trajectory optimization algorithm. Motivated by this observation, we tackle the initialization problem for trajectory optimization via policy optimization. To optimize a policy, we propose a guided policy search method that has two key components: i) Trajectory update; ii) Policy update. The trajectory update involves offline solutions of a large number of trajectory optimization problems from different initial states via Sequential Convex Programming (SCP). Here we take a single SCP step to generate the trajectory iterate for each problem. In conjunction with these iterates, we also generate additional trajectories around each iterate via a feedback control law. Then all these trajectories are used by a stochastic gradient descent algorithm to update the neural network policy, i.e., the policy update step. As a result, the trained policy makes it possible to generate trajectory candidates that are close to the optimality and feasibility and that provide excellent initial guesses for the trajectory optimization methods. We validate the proposed method via a real-world 6-degree-of-freedom powered descent guidance problem for a reusable rocket.
\end{abstract}


\section{Introduction}

Trajectory optimization methods have been developed to solve nonconvex trajectory generation problems, e.g., optimal control problems with nonlinear dynamics and nonconvex constraints. Many different algorithms for trajectory optimization have been developed, such as sequential convex programming (SCP) \cite{szmuk2018successive}, differential dynamics programming (DDP) \cite{koenemann2015whole}, and Pontryagin's maximum principle (PMP) \cite{kim2002nonlinear}. With these developed methods, numerous studies have applied trajectory optimization to various applications including control of mobile robots \cite{kim2017path}, autonomous robotic manipulation \cite{kumar2016optimal}, and powered descent guidance for spacecraft and reusable rockets \cite{szmuk2018successive}.  

A key input of many trajectory optimization algorithms including SCP and DDP is an initial trajectory guess, which is used to start the iterative process \cite{malyuta2021convex}. This initial trajectory guess is observed to impact the number of iterations for the optimization to converge, and hence it can be critical for real-time applications. Thus, designing a reliable and efficient initial trajectory guess can be beneficial for the fast convergence of trajectory optimization. The prior methods for the initial guess, such as a straight-line interpolation, however, are based on engineering heuristics motivated by the user's experiences. For a more systematic method, we propose a policy optimization-based machine learning approach for trajectory optimization initialization. The policy optimization aims to acquire a neural network policy that can give the optimal control action as a function of the  state. If we can obtain a proper neural network policy through policy optimization, we can then generate trajectory candidates that serve as excellent initial guesses for the trajectory optimization methods by integrating the system dynamics with control inputs from the optimal policy.

One way to optimize and acquire the  policy is to employ model-free
reinforcement learning (RL). In recent years, much attention has been drawn to the
model-free RL because of its applicability and performance
for various complex problems. For example, \cite{furfaro2020adaptive}
deployed the actor-critic based policy gradient RL method for a
zero-effort-miss/zero-effort-velocity feedback guidance problem. The
proximal policy optimization RL method was used for the satellite
rendezvous missions in \cite{broida2019spacecraft}. One disadvantage
of applying model-free RL is sample inefficiency. As model-free
RL methods do not utilize the system models, many data samples are
required to train a policy. Another disadvantage is that
it is still not straightforward to impose the state and input constraints necessary for safety-critical systems like space vehicles
and self-driving cars.

Another way to optimize a policy is to use imitation learning \cite{osa2018algorithmic} in which trajectory
optimization methods play a role in giving samples, and the neural net
policy is trained to these samples through supervised learning.
The PMP method was used by \cite{sanchez2018real} to generate trajectory samples to train a neural network policy via supervised learning. This method was validated on the soft-landing problem with a 3-degree-of-freedom (DoF) rocket.
Real-time trajectory generation for asteroid landings was conducted
by imitation learning in \cite{cheng2020real}. 
The downside of imitation learning
is that the errors in supervised learning can be compounded along
the trajectory over a time horizon. Then, this accumulated error can
result in poor propagated trajectories over a long time horizon especially for
higher-dimensional systems.

Motivated by the shortcoming of imitation learning, interactive
optimization approaches, also known as guided policy search, have recently been studied in \cite{levine2016end,mordatch2015interactive}.
These methods formulate the entire policy optimization problem as two separated optimizations: one is to obtain solutions from different initial conditions with the trajectory optimization routine, and the other is to train the neural net policy using supervised learning with sample trajectories generated in the previous step. Unlike imitation learning, the trajectory optimization in the first step has an additional constraint that enforces the updated trajectory to be close enough to the trajectory from the policy. Imposing this constraint can help the single neural network policy to reproduce all trajectories with the long-time horizon performance. 

In this paper, we propose a guided policy search algorithm that leverages convex optimization. We repeatedly solve two optimization sub-problems: the first step aims to solve a convex approximation of the nonlinear optimal control problem with an additional penalty on control input deviation from the policy. In this step, we formulate the sub-problem as a convex optimization problem, which is motivated by the sub-problem in a penalized trust region (PTR) algorithm which is one of SCP methods. Then, we generate additional samples around the solution of the convex sub-problem using a feedback gain obtained by a time-varying linear quadratic regulator (LQR). The second step is to train the neural net policy using supervised learning with the generated samples. Then, the two optimizations with the sample generation are repeatedly solved until the policy converges. We validate our method on the minimum-fuel powered descent guidance for a reusable rocket. The simulation results show that the policy can generate trajectory candidates which lead to the fast and reliable convergence of the PTR method.

\subsection{Related works}

The scheme of the proposed method is related to the previous research
in the guided policy search methods, but there are technical differences.
The works of \cite{levine2016end,carius2020mpc} utilized the DDP method and \cite{mordatch2015interactive}
worked with Newton's method.
On the other hand, we employ a convex approximation of the original
nonconvex problem instead of solving the nonconvex optimal control problem in each iteration.
There are two advantages of this approach. First, solving the convex
approximation instead of solving the nonconvex trajectory optimization problem is more computationally tractable. Since this trajectory optimization
is repeated as the first step in each iteration, solving the nonconvex
problem is too costly. More importantly, leveraging convex optimization
makes it straightforward to impose state and input constraints while
most earlier works did not impose the constraints. For instance, the work in \cite{kim2018vision} imposed the state constraints by using an augmented Lagrangian formulation, but the constraint exists for the training phase, and the final policy does not activate the constraint. The work in the imitation learning \cite{carius2020mpc} considered both inequality and equality constraints, but the inequality constraint are implemented as a penalty via a barrier function, not as a hard constraint.

\subsection{Contributions}

We can summarize our main contributions as follows:
\begin{itemize}
\item The proposed method makes it straightforward to impose state and input constraints in the framework of guided policy search by leveraging convex optimization.
\item The trained policy by the proposed method can provide a better initial guess in comparison to the conventional initial guesses used for trajectory optimization.
\item The proposed approach is validated via minimum-fuel 6-DoF powered descent guidance problem for reusable rockets, which is a real-world aerospace application.
\end{itemize}

\subsection{Outlines}

The rest of this paper is organized as follows: We present the problem formulation and the proposed method in section \ref{sec2}. In section \ref{sec3}, we describe the details of the 6-DoF powered descent guidance for a reusable rocket. In section \ref{sec4}, we show numerical results with our proposed method. Concluding remarks are provided in \ref{sec5}.

\section{Policy Optimization via Penalized Trust Region}
\label{sec2}

\subsection{Problem formulation}

In this paper, we consider a deterministic continuous-time policy optimization problem of the following form:
\begin{align}
\min_{\genfrac{}{}{0pt}{1}{\theta,x^{i}(\cdot),u^{i}(\cdot),}{i=1,\ldots,N}}\: & \sum_{i=1}^{N}\int_{0}^{t_{f}}J(t,x^{i}(t),u^{i}(t))\mathrm{d}t\label{eq:cost}\\
\text{s.t.}\qquad
 & \dot{x}^{i}(t)=f(t,x^{i}(t),u^{i}(t)),\label{eq:model}\\
 & x^{i}(t)\in\mathcal{X}(t),\quad u^{i}(t)\in\mathcal{U}(t),\label{eq:convexconst}\\
 & s(t,x^{i}(t),u^{i}(t))\leq0,\label{eq:nonconvexconst}\\
 & \,\quad \qquad \qquad \qquad \qquad \forall t\in[0,t_f],\nonumber\\
 & u^{i}(t)=\pi_{\theta}(x^{i}(t_{k})),\quad\forall t\in[t_{k},t_{k+1}),\label{eq:policyconst}\\
 & \quad\quad\quad\quad\quad\quad\quad\quad\quad  k=0,\ldots,K-1, \nonumber \\
 & x^{i}(0)=x_{\text{{init}}}^{i},\quad x^{i}(t_{f})=x_{\text{final}}^{i},\label{eq:statefinal}
\end{align}
where the subscript $i$ indexes each trajectory starting from the initial conditions $x_{\text{{init}}}^{i}$ and ending at the final state $x_{\text{final}}^{i}$. The user-specified parameter $N$ is the number of trajectories and $t_{f}$ is the fixed final time. The vector $x(\cdot)\in\mathbb{R}^{n_{x}}$
is the state and $u(\cdot)\in\mathbb{R}^{n_{u}}$ is the control input.
The function $f:\mathbb{R}\times\mathbb{R}^{n_{x}}\times\mathbb{R}^{n_{u}}\rightarrow\mathbb{R}^{n_{x}}$
represents the (nonlinear) system dynamics. The constraints for the
optimization are given in \eqref{eq:convexconst} and \eqref{eq:nonconvexconst}.
The convex sets $\mathcal{X}(t)$ and $\mathcal{U}(t)$ are state
constraints and input constraint sets, respectively. We define a function
$s:\mathbb{R}\times\mathbb{R}^{n_{x}}\times\mathbb{R}^{n_{u}}\rightarrow\mathbb{R}^{n_{s}}$
to illustrate nonconvex constraints. The function $\pi_{\theta}:\mathbb{R}^{n_{x}}\rightarrow\mathbb{R}^{n_{u}}$
represents a policy parameterized by a neural network, and
$\theta$ is the vector of weights for the neural network. We assume that the
input from the policy is piecewise  constant between the sampling times as in \eqref{eq:spacedtemporalgrid}, which is the form of zeroth order hold \cite{malyuta2019discretization}. 
\begin{align}
t_{k} & = \frac{k}{K-1}t_{f},\quad k=0,\ldots,K-1\label{eq:spacedtemporalgrid}.
\end{align}
The policy
optimization aims to optimize the neural net policy that approximately reproduce trajectories produced by trajectory optimization subject to different initial conditions. After optimizing the policy, we can obtain
a new trajectory with the trained policy by propagating the following system
dynamics
\begin{align*}
u(t) & =\pi_{\theta}(x(t_{k}))\quad\forall t\in[t_{k},t_{k+1}),\\
\dot x(t) &= f(t,x(t),u(t)).
\end{align*}

\subsection{Penalized trust region}

Before diving into the proposed method, we review the PTR method which
is a key component in this research. The PTR method, which is a type of SCP, aims to solve (nonconvex)
continuous-time optimal control problems that can be formulated as
follows:
\begin{align*}
\min_{\genfrac{}{}{0pt}{1}{x^{i}(\cdot),u^{i}(\cdot),}{i=1,\ldots,N}}\: & \sum_{i=1}^{N}\int_{0}^{t_{f}}J(t,x^{i}(t),u^{i}(t))\mathrm{d}t\\
\text{s.t.}\quad & \eqref{eq:model}, \eqref{eq:convexconst}, \eqref{eq:nonconvexconst}, \eqref{eq:statefinal}.
\end{align*}

To solve the above optimization, the PTR iteratively solves a convex sub-problem.
Every sub-problem's solution provides a new reference trajectory to the algorithm. After linearizing and discretizing around this reference,
the convex sub-problem can be formulated as follows:
\begin{align}
\min_{\substack{ x_{k}^{i},u_{k}^{i},\nu_{k}^{i}, \\ i=1,\ldots,N, \\ k=0,\ldots,K-1}}\: & \sum_{i=1}^{N}\sum_{k=0}^{K-1}J(t_{k},x_{k}^{i},u_{k}^{i})+J_{vc}(\nu_{k}^{i})+J_{tr}(x_{k}^{i},u_{k}^{i})\label{eq:SOCP_cost}\\
\text{s.t.}\quad\, & x_{k+1}^{i}=A_{k}^{i}x_{k}^{i}+B_{k}^{i}u_{k}^{i}+z_{k}^{i}+\nu_{k}^{i},\label{eq:SOCP_model}\\
 & x_{k}^{i}\in\mathcal{X}(t_k),\quad u_{k}^{i}\in\mathcal{U}(t_k),\label{eq:SOCP_convex}\\
 & C_{k}^{i}x_{k}^{i}+D_{k}^{i}u_{k}^{i}+r_{k}^{i}\leq0,\label{eq:SOCP_nonconvex}\\
& x_{0}^{i}=x_{\text{{init}}}^{i},\quad x_{K-1}^{i}=x_{\text{final}}^{i}\label{eq:SOCP_ini}, 
\end{align}
where \eqref{eq:SOCP_nonconvex} is the first-order approximation
of the original nonconvex constraint in \eqref{eq:nonconvexconst},
and $A_{k}^{i},B_{k}^{i},z_{k}^{i}$ represent the linearized model
of the nonlinear dynamics given in \eqref{eq:model} around a reference
trajectory denoted by $\bar{x}_{k}^{i},\bar{u}_{k}^{i}$.


Besides the original cost $J$, we have
two additional penalties in (\ref{eq:SOCP_cost}): virtual control and trust region \cite{malyuta2021convex}.
Virtual control serves to recover from cases where the linearized sub-problem becomes infeasible.
Trust regions maintain the optimization in the vicinity of the reference trajectory, where the linearization is most accurate. They can
be formulated as follows:
\begin{align}
J_{vc}(\nu_{k}^{i}) & =w_{\nu}\norm{\nu_{k}^{i}},\nonumber \\
J_{tr}(x_{k}^{i},u_{k}^{i}) & =w_{tr}(\norm{x_{k}^{i}-\bar{x}_{k}^{i}}_{2}^{2}+\norm{u_{k}^{i}-\bar{u}_{k}^{i}}_{2}^{2}),\label{eq:ptrtr}
\end{align}
where $\nu_{k}^{i}$ is the virtual control variable. Then, the PTR algorithm terminates when the following two stopping criteria are satisfied: 
\begin{align}
\sum_{k=0}^{K-1}\norm{\nu_{k}}_{1}\leq\epsilon_{\nu},\quad\sum_{k=0}^{K-1}\norm{x_{k}^{i}-\bar{x}_{k}^{i}}_{2}^{2}+\norm{u_{k}^{i}-\bar{u}_{k}^{i}}_{2}^{2}\leq\epsilon_{tr},\label{eq:PTR_termination}
\end{align}
where $\epsilon_{\nu}$ and $\epsilon_{tr}$ are user-specified tolerance parameters that are typically small numbers. More details about the PTR algorithm can be found in \cite{szmuk2017successive,szmuk2018successive}. 

\subsection{Guided policy search via SCP}

\label{subsec:PTR}
\begin{figure}
\begin{center}
\includegraphics[width=8.3cm]{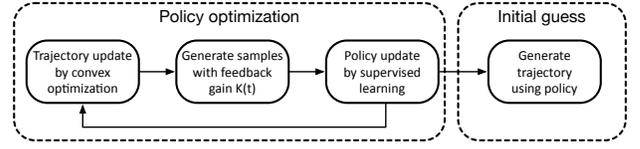}  
\caption{
The block diagram of the proposed method. The method repeatedly solves two optimizations for the update of trajectory and policy. Between two optimizations, there is a step for generating samples that will be used for the policy update. Then, the trained policy generates a trajectory that will be an initial guess for trajectory optimization.
} 
\label{fig:blockdiagram}
\end{center}
\end{figure}

The policy optimization given in \eqref{eq:cost}-\eqref{eq:statefinal}
is solved in the proposed guided policy search algorithm. To optimize the policy, we repeatedly solve two optimization steps: a trajectory update step and a policy
update step. In the trajectory update step, we solve the convex sub-problem
that approximates the original nonconvex problem. Next, we generate trajectory samples using the solution of the sub-problem and feedback gain matrix that we obtain using a linear time-varying LQR.
With the trajectory samples, we train a neural net policy
by supervised learning. The trajectory and policy update steps are repeated in sequence until the policy converges. The block diagram for the proposed approach is given in Fig.~\ref{fig:blockdiagram}.

The convex sub-problem for the trajectory update step is formulated as follows:
\begin{align}
\min_{\substack{ x_{k}^{i},u_{k}^{i},\nu_{k}^{i}, \\ i=1,\ldots,N, \\ k=0,\ldots,K-1}}\: & \sum_{i=1}^{N}\sum_{k=0}^{K-1}J(t_{k},x_{k}^{i},u_{k}^{i})+J_{vc}(\nu)+J_{trp}(u_{k}^{i})\nonumber \\
\text{s.t.}\quad\: & \eqref{eq:SOCP_model}-\eqref{eq:SOCP_ini}.\label{eq:subproblemPO}
\end{align}
The formulation of the sub-problem is similar to
that of the PTR given in \eqref{eq:SOCP_cost}-\eqref{eq:SOCP_ini}.
The only difference is a penalty term $J_{trp}(u)$ for a trust region
in \eqref{eq:trp}.

\begin{equation}
J_{trp}(u)=w_{trp}\norm{u_{k}^{i}-\pi_{\theta}(\bar{x}_{k}^{i})}_2^2.\label{eq:trp}
\end{equation}

While the trust region  \eqref{eq:ptrtr} in the PTR is for the difference
with the reference trajectory, the new trust region in \eqref{eq:trp}
represents the deviation of the control input from the output of the policy.
This trust region with the policy ensures that the solution from the
convex sub-problem does not deviate too far from the trajectory
made by the policy $\pi_{\theta}$. After solving
the trajectory update step, we can get the updated trajectories as $\{\{\hat{x}_{k}^{i},\hat{u}_{k}^{i},\hat{\nu}_{k}^{i}\}_{k=0}^{K-1}\}_{i=1}^{N}$.


The next step is to generate sample trajectories using the trajectory resulting from the trajectory
update step. In \cite{levine2016end}, the authors suggested the way
of generating samples around the given trajectory using a time-varying
feedback gain $K(k)$. By giving additional samples, we can increase the amount of data
and make the neural network be more robust to the small perturbation in the trajectory. Following a similar approach, first, we obtain the feedback gain by a linear time-varying LQR solution with the linearized model
given in \eqref{eq:SOCP_model}. Then, we choose $S$ new initial
states $x_{\text{init}}^{i,s}$ around the original ones $x_{\text{init}}^{i}$
by adding Gaussian noise, and generate $S$ new trajectories with the
linearized model. This can be formulated as follows:
\begin{align}
x_{\text{init}}^{i,s} & =x_{\text{init}}^{i}+\epsilon^{i,s},\epsilon^{i,s}\sim\mathcal{N}(0,\sigma_{i}^{2}),\nonumber \\
u_{k}^{i,s} & =\hat{u}_{k}^{i}+K^{i}(k)(x_{k}^{i,s}-\hat{x}_{k}^{i}),\nonumber \\
x_{k+1}^{i,s} & =A_{k}^{i}x_{k}^{i,s}+B_{k}^{i}u_{k}^{i,s}+z_{k}^{i}+\hat{\nu}_{k}^{i}.\label{eq:generatesamples}
\end{align}
If the states include quaternions, the noise is added to roll, pitch, and yaw angles that are converted from the quaternions, and then these angles are transformed again into the quaternions. 

One of the disadvantages of the above approach is that the linearized trajectory samples might violate the original constraints of the sub-problem. To mitigate this issue, we gradually decrease
the magnitude of variance $\sigma_i^{2}$ in each iteration to minimize the constraint violation. A promising future direction to handle the constraint violation is to employ funnel synthesis approaches in which the feedback gain is calculated with the consideration of the constraints \cite{reynolds2021funnel}. 

Finally, we apply supervised learning with the generated samples to
train the policy represented by a neural network in the policy
update step. We optimize the objective function given in \eqref{eq:policyobj}
using stochastic gradient descent (SGD) over state and input pairs in the trajectory samples.
\begin{align}
J_{p}(\theta)=\sum_{s=1}^{S}\sum_{i=1}^{N}\sum_{k=0}^{K-1}\norm{u_{k,s}^{i}-\pi_{\theta}(x_{k,s}^{i})}\label{eq:policyobj}.
\end{align}

The overall approach is summarized in Algorithm \ref{alg:alg1}:

\begin{algorithm}
\begin{algorithmic}[1]
\State{initialize parameters $\theta$ of network for policy}
\While{not converged}
	\State{set the reference trajectory $\bar{x}_{k}^{i},\bar{u}_{k}^{i}$}
	\State{linearize and discretize model and constraints}
	\State{update trajectories $\hat{x}_{k}^{i},\hat{u}_{k}^{i},\hat{\nu}_{k}^{i}$ \eqref{eq:subproblemPO}}
	\State{generate samples $x_{k,s}^{i},u_{k,s}^{i}$ \eqref{eq:generatesamples}}
	\State{optimize parameters of the policy $\theta$ \eqref{eq:policyobj}}
\EndWhile
\end{algorithmic}
\caption{Policy optimization via PTR}
\label{alg:alg1}
\end{algorithm}

The parameters $\theta$ of the neural net policy are initialized by the Xavier initialization method proposed in \cite{glorot2010understanding}. In the first iteration, the reference trajectory is obtained by the straight-line interpolation \cite{malyuta2021convex}. Since having the penalty with the randomly initialized neural net policy in the first iteration can slow down the entire process, the convex sub-problem in \eqref{eq:subproblemPO} uses the trust region penalty with the input $\bar{u}$ from the reference trajectory instead of the input from the policy $\pi_{\theta}(\bar{x})$ in \eqref{eq:trp}. Then, the convex sub-problem's solution provides a new reference trajectory for the next iteration, and the trust region penalty with the policy is employed from the second iteration. 

The proposed algorithm determines its stopping criteria with validation data. The validation data is the set of trajectories generated from new initial states that are not identical with the states $x_{\text{init}}^i$ used in the training. After training the policy in every iteration, we evaluate the average cost $J$ of trajectories starting from these new initial states in the validation data. If the difference between the average cost of the current iteration and that of the last iteration is less than a user-specified tolerance  $\epsilon_J$, the method terminates the process. Note that the trajectories in the validation set are used purely for evaluating the policy performance and not for the learning phase.

\section{Application to Minimum-fuel Powered Descent Guidance}
\label{sec3}
The method is validated on the problem of powered descent guidance for a reusable rocket. This section outlines the problem formulation. The readers can find more details in \cite{szmuk2018successive}.

\subsection{Dynamics and kinematics}

The system dynamics for a 6-DoF rocket in an $X$(East)-$Y$(North)-$Z$(Up) inertial coordinate frame is given by:
\begin{align}
\dot{m}(t) & =-\alpha_{\dot{m}}\norm{T_{\mathcal{B}}(t)}_{2},\label{eq:const_pgdmass}\\
\dot{r}_{\mathcal{I}}(t) & =v_{\mathcal{I}}(t),\\
\dot{v}_{\mathcal{I}}(t) & =\frac{1}{m(t)}C_{\mathcal{I}/\mathcal{B}}(t)T_{\mathcal{B}}(t)+g_{\mathcal{I}},\\
\dot{q}_{\mathcal{B}/\mathcal{I}}(t) & =\frac{1}{2}\Omega(\omega_{\mathcal{B}}(t))q_{\mathcal{B}/\mathcal{I}}(t),\\
J_{\mathcal{B}}\dot{\omega}_{\mathcal{B}}(t) & =[r_{T,\mathcal{B}}\times]T_{\mathcal{B}}(t)-[\omega_{\mathcal{B}}(t)\times]J_{\mathcal{B}}\omega_{\mathcal{B}}(t),
\end{align}
where $m(\cdot)\in\mathbb{R}_{++}$ is the mass of the vehicle, $\alpha_{\dot{m}}$
is the proportionality constant, and $T_{\mathcal{B}}(\cdot)\in\mathbb{R}^{3}$
is the thrust vector expressed in the body-fixed frame denoted by
$\mathcal{F}_{\mathcal{B}}$. The position, velocity, and
constant gravitational acceleration vector are given by $r_{\mathcal{I}}(\cdot)\in\mathbb{R}^{3}$,
$v_{\mathcal{I}}(\cdot)\in\mathbb{R}^{3}$, and $g_{\mathcal{I}}\in\mathbb{R}^{3}$,
respectively. The quaternions are represented by $q_{\mathcal{B}/\mathcal{I}}\in\mathcal{S}^{3}$ to parameterize the attitude of $\mathcal{F}_{\mathcal{B}}$
relative to a inertial reference frame $\mathcal{F}_{I}$, and $C_{\mathcal{I}/\mathcal{B}}(t)\in SO(3)$
is the direction cosine matrix to show the attitude transformation
from $\mathcal{F_{B}}$ to $\mathcal{F}_{\mathcal{I}}$. The angular
velocity vector is denoted by $\omega_{\mathcal{B}}(t)\in\mathbb{R}^{3}$,
and $J_{\mathcal{B}}\in\mathbb{S}^{3}$ is the moment of inertia.
The notations given as $[\xi\times]\in\mathbb{R}^{3\times 3}$ and $\Omega(\xi)\in\mathbb{R}^{4\times 4}$
represent skew-symmetric matrices for some $\xi\in\mathbb{R}^{3}$.

\subsection{State and input constraints}

The rocket is subject to the following state constraints:
\begin{align}
m_{dry} & \leq m(t),\label{eq:const_start}\\
\norm{\omega_{\mathcal{B}}(t)}_{2} & \leq\omega_{max},\\
\tan\gamma_{gs}\norm{H_{12}^{T}r_{\mathcal{I}}(t)}_{2} & \leq\mathbf{e}_{3}^{\top}r_{\mathcal{I}}(t),\\
\cos\theta_{max} & \leq1-2(q_{2}^{2}(t)+q_{3}^{2}(t)),
\end{align}
where $m_{dry}$ is the dry mass and $\omega_{max}$ is the maximum
angular rate. The maximum glide slope angle is $\gamma_{gs}$ and
$H_{12}$ is the indicator for first and second elements that is
written as $H_{12}=[\mathbf{e}_{1},\mathbf{e}_{2}]$ where we use $\mathbf{e}_{j}$
to illustrate the unit vector for the $j$-th axis. The maximum tilt angle of the vehicle is $\theta_{max}$, and $q_{2}$ and $q_{3}$
are second and third components in the quaternion $q_{\mathcal{B}/\mathcal{I}}$.

The input constraints for the vehicle are given by:
\begin{align}
0<T_{min}\leq\norm{T_{\mathcal{B}}(t)}_{2} & \leq T_{max},\\
\cos\delta_{max}\norm{T_{\mathcal{B}}(t)}_{2} & \leq\mathbf{e}_{3}^{\top}T_{\mathcal{B}}(t),\label{eq:const_pgdgimbalangle} 
\end{align}
where $T_{min}$ and $T_{max}$ are the minimum and maximum thrust,
respectively. The maximum gimbal angle is $\delta_{max}$.
\subsection{Boundary conditions}

The state values are subject to initial and terminal boundary constraints, given as follows:
\begin{align}
m(0) &= m_{wet},\\
r_{\mathcal{I}}(0) &= r_{\mathcal{I},i},\quad r_{\mathcal{I}}(t_{f})=r_{\mathcal{I},f},\\
v_{\mathcal{I}}(0) & =v_{\mathcal{I},i},\quad v_{\mathcal{I}}(t_{f})=v_{\mathcal{I},f},\\
q_{\mathcal{B}/\mathcal{I}}(0) & =q_{\mathcal{B}/\mathcal{I},i},\quad q_{\mathcal{B}/\mathcal{I}}(t_{f})=q_{\mathcal{B}/\mathcal{I},f},\\
\omega_{\mathcal{B}}(0) & =\text{\ensuremath{\omega_{\mathcal{B},i}}},\quad\omega_{\mathcal{B}}(t_{f})=\text{\ensuremath{\omega_{\mathcal{B},f}}}.\label{eq:const_pgdlastboundarycondtion}
\end{align}

\subsection{Problem statement}

With the derived constraints, the minimum-fuel optimal control problem
for the 6-DoF rocket soft-landing can be formulated as follows:
\begin{align*}
\min_{T_{\mathcal{B}}(t)}\: & -m(t_{f})\\
\text{s.t.}\: & \eqref{eq:const_pgdmass}-\eqref{eq:const_pgdlastboundarycondtion}.
\end{align*}

\section{Simulation Results}
\label{sec4}
\begin{figure*}
\begin{center}
\includegraphics[width=16cm]{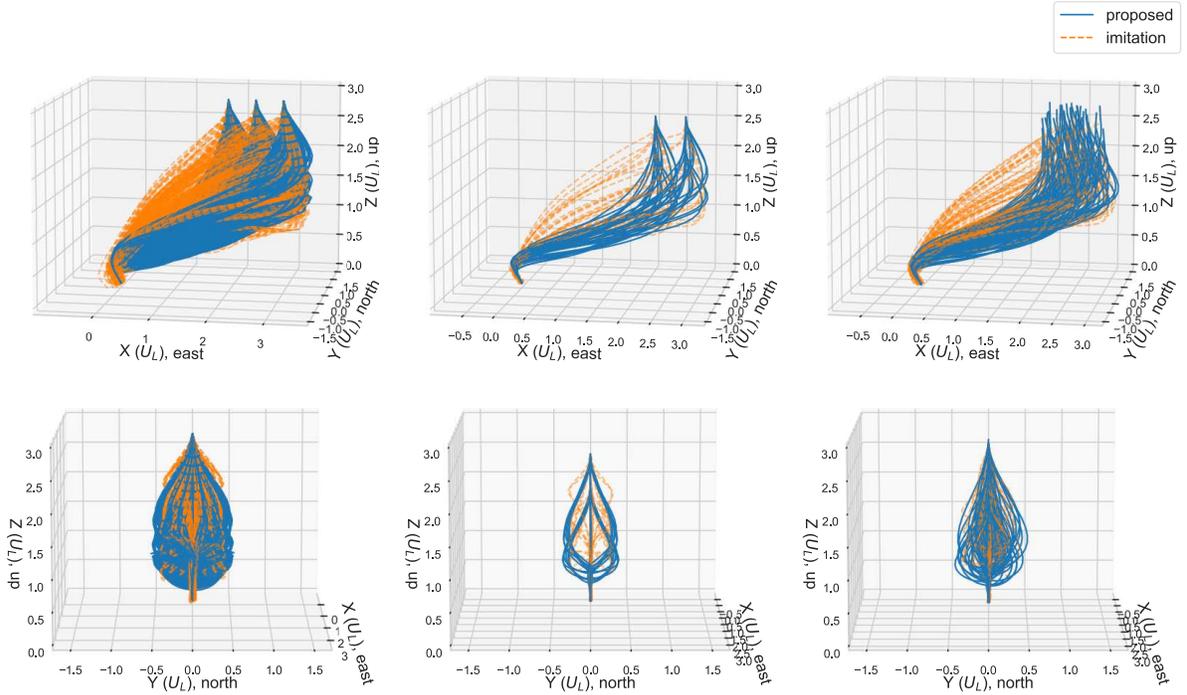}  
\caption{
The trajectory comparison between the proposed approach and the imitation learning method. The figures in the first column represent the trajectories starting from the initial states in the training set. The trajectories in the second column images and in the third column images start from the initial states in the validation set and the test set, respectively. The first and second rows show the same set of trajectories viewed from different angles. It is visible that the proposed approach generates less infeasible trajectories than the imitation learning method in terms of the boundary condition for the position.
} 
\label{fig:trajectory_comparison}
\end{center}
\end{figure*}
This section presents numerical results for the powered descent guidance problem. The simulations aim to answer the following questions:
\begin{itemize}
\item Can the trained policy generate trajectories from different initial conditions with long time horizon performance?
\item Can a trajectory optimization algorithm converge faster with the initial guess provided by the trained policy instead of the conventional initialization method?
\end{itemize}

Initial states for the simulation are selected by the following manner: the first and third components of the position $r_\mathcal{I}$ are allowed to take 3 values as $\{2 (U_L),2.5 (U_L),3 (U_L)\}$, the first and second elements of the velocity $v_\mathcal{I}$ take a value among $\{-0.1 (U_L/U_T),0.0 (U_L/U_T),0.1 (U_L/U_T)\})$, and roll and pitch angles take 3 values as $\{-15 (\text{\textdegree}),0 (\text{\textdegree}),15 (\text{\textdegree})\}$. The roll and pitch angles are then transformed to quaternions. The notations $U_T$,$U_L$, and $U_M$ are nondimensionalized time, length, and mass units. We fix the initial mass $m$ as $2.0~(U_M)$, and the third component of the velocity as $-1.0~(U_L/U_T)$, and the remaining states as zero. Among 729 ($3^6$) combinations for the initial states, we exclude 22 initial states with which the PTR does not converge, and then a total of 707 initial states are selected to generate the training set. Each initial state is a starting point of a trajectory, so the number of trajectories $N$ in \eqref{eq:cost}-\eqref{eq:statefinal} is 707. 

We also select 36 and 100 new initial states to generate a validation set and a test set, respectively. In the initial states in validation set, the first and third components of the position $r_\mathcal{I}$ take 2 values as $\{2.25 (U_L),2.75 (U_L)\}$, and roll and pitch angles take a value within $\{-15 (\text{\textdegree}),0 (\text{\textdegree}),15 (\text{\textdegree})\}$. The initial mass $m$, the third component of the velocity, and other states are fixed as $2.0~(U_M)$, $-1.0~(U_L/U_T)$, and zero, respectively. The initial states in the test set follow a uniform random distribution inside the rectangular region defined by the ranges of initial state values used for the training set. We illustrate the initial states and trajectories by the trained policy for the training, validation, and test sets in Fig.~\ref{fig:trajectory_comparison}. 

The results in this section use a neural network with three hidden layers and 50 nodes in each layer to represent the policy. Each layer except the last one has the rectified linear unit layer as an activation function. The batch size of the training is 512. Other parameters for the simulation are given in Table \ref{tb:table_example}. The proposed approach reaches its stopping criteria at 5 iterations. 
\begin{table}
\caption{Simulation Parameters}
\label{tb:table_example}
\begin{center}
\begin{tabular}{cccc}
\toprule 
Parameter & Value (Unit) & Parameter & Value (Unit)\tabularnewline
\midrule 
$g_{\mathcal{I}}$ & [0,0,-1] ($U_{L}/U_{T}^{2}$) & $\delta_{max}$ & 20 ($\text{\textdegree}$)\tabularnewline
\midrule 
$J_{xx}$ & 0.186 ($U_{M}U_{L}^{2}$) & $\theta_{max}$ & 90 ($\text{\textdegree}$)\tabularnewline
\midrule 
$J_{yy}$ & 0.186 ($U_{M}U_{L}^{2}$) & $\omega_{max}$ & 60 ($\text{\textdegree}/U_{T}$)\tabularnewline
\midrule 
$J_{zz}$ & 0.00372 ($U_{M}U_{L}^{2}$) & $\gamma_{gs}$ & 20 ($\text{\textdegree}$)\tabularnewline
\midrule 
$J_{\mathcal{B}}$ & 
\begin{tabular}{@{}c@{}}$\text{diag}([J_{xx},J_{yy},J_{zz}])$ \\ ($U_{M}U_{L}^{2}$) \end{tabular}
& $r_{T,\mathcal{B}}$ & 0.25 ($U_{L})$\tabularnewline
\midrule 
$\alpha_{\dot{m}}$ & 0.01 ($U_{T}/U_{L}$) & $w_{trp}$ & 10\tabularnewline
\midrule 
$m_{wet}$ & 2.0 ($U_{M}$) & $w_{\nu}$ & $10^{4}$\tabularnewline
\midrule 
$m_{dry}$ & 1.0 ($U_{M}$) & $N$ & 707\tabularnewline
\midrule 
$T_{max}$ & 6.0 ($U_{M}U_{L}/U_{T}^{2}$) & $K$ & 31\tabularnewline
\midrule 
$T_{min}$ & 1.5 ($U_{M}U_{L}/U_{T}^{2}$) & $S$ & 20\tabularnewline
\midrule 
$\epsilon_J$ & $10^{-4}$ & $t_{f}$ & 5 ($U_{T}$)\tabularnewline
\midrule 
$r_{\mathcal{I},f}$ 
& 
\begin{tabular}{@{}c@{}}[0,0,0] \\  ($U_L$) \end{tabular}
& 
$v_{\mathcal{I},f}$ 
& 
\begin{tabular}{@{}c@{}} [0,0,-0.1] \\ ($U_L/U_T$)\end{tabular}\tabularnewline
\midrule 
$q_{\mathcal{B},f}$ & [1,0,0,0] (-) & $\omega_{\mathcal{B},f}$ & [0,0,0] (\textdegree$/U_T$)\tabularnewline
\bottomrule
\end{tabular}
\end{center}
\end{table}

\subsection{Comparison to imitation learning}

The comparisons between our approach and an imitation learning method in terms of cost, constraints, and boundary condition of the trajectories are given in here. The latter approach decouples the PTR trajectory generation phase from the neural network training. Effectively, a standard PTR method generates complete trajectories from the initial states in the training and validation sets. Around each trajectory generated by the PTR method with the training set, we generate 100 samples by following the same procedure with our approach, which is given in \eqref{eq:generatesamples}. In the proposed method, 20 samples are generated, and the number of iterations is 5, so the number of samples around each trajectory is 100, which is equal to that of the imitation learning method. Subsequently, supervised learning is used to independently train a neural network on the samples. In every epoch in training, we evaluate the policy with the objective given in (\ref{eq:policyobj}) using the state and input pairs in the trajectories starting from the initial state in the validation set. Then, the network in the epoch that shows the minimum objective value is selected as the final policy for the imitation learning method. We employ the same network structure used in the proposed approach for a fair comparison. 

Table~\ref{tb:details} compares the two algorithms in terms of cost, constraints, and boundary conditions. The first row represents the average cost of all trajectories. The average maximum constraint violation of trajectories is given in the second row where $c_k \in \mathbb{R}^7$ is a normalized constraint violation defined as
\begin{align}
\label{eq:normalized_const}
c_k = 
\begin{bmatrix}
(m_{dry}-m_k)^{+} / m_{dry} \\
(\omega_k-\omega_{max})^{+} / \omega_{max} \\
(\gamma_{gs}-\gamma_k)^{+} / \gamma_{gs} \\
(\theta_k-\theta_{max})^{+} / \theta_{max} \\
(T_{min}-\norm{T_{k,\mathcal{B}}})^{+} / T_{min} \\
(\norm{T_{k,\mathcal{B}}}-T_{max})^{+} / T_{max} \\
(\delta_k-\delta_{max})^{+} / \delta_{max}
\end{bmatrix},
k=0,\ldots,K-1. \nonumber
\end{align}
The function $(x)^+=\max (0,x)$ is an unit ramp function. The last 4 rows show the average boundary condition violations for the position, velocity, attitude, and angular velocity separately. The quaternions are transformed into the Euler angle $\theta_e$, and then the violation of attitude is calculated. The results show that the generated trajectories by our approach are less infeasible in terms of the satisfaction of both constraints and boundary conditions. Even though the imitation learning method has less cost than that of our approach, it is expected that a trajectory that does not respect the problem's constraints will achieve a lower cost. We describe the trajectories generated by both our proposed approach and the imitation learning in Fig.~\ref{fig:trajectory_comparison}. We can see that the proposed approach generates less infeasible trajectories in terms of the position boundary condition than the imitation learning method.

\begin{table}
\caption{Evaluation of Cost, Constraints and Boundary Condition}
\label{tb:details}
\begin{tabular}{ccccc}
\toprule 
 & \multicolumn{2}{c}{Validation} & \multicolumn{2}{c}{Test}\tabularnewline
\midrule 
Property (Unit) & Proposed  & Imitation & Proposed & Imitation\tabularnewline
\midrule 
\begin{tabular}{@{}c@{}} $-m(t_{f})$ \\ ($U_{M}$) \end{tabular} & -1.8591 & -1.8648 & -1.8610 & -1.8670\tabularnewline
\midrule
\midrule 
\begin{tabular}{@{}c@{}}  $\max_k\norm{c_k}_\infty$ \\ (-) \end{tabular} 
& 0.0031 & 0.0615 & 0.0361 & 0.2112 \tabularnewline
\midrule
\midrule 
\begin{tabular}{@{}c@{}} $\norm{r_{\mathcal{I}}(t_{f})-r_{\mathcal{I},f}}_2$
\\ ($U_L$) \end{tabular} & 0.0050 & 0.0240 & 0.0081 & 0.0315
\tabularnewline
\midrule
\begin{tabular}{@{}c@{}}
$\norm{v_{\mathcal{I}}(t_{f})-v_{\mathcal{I},f}}_2$ 
\\ ($U_L/U_T$) \end{tabular}
& 0.0059 & 0.0567 & 0.0102 & 0.0440
\tabularnewline
\midrule
\begin{tabular}{@{}c@{}}
$\norm{\theta_e(t_{f})-\theta_{e,f}}_2$
\\ ($\text{\textdegree}$) \end{tabular}
& 0.3074 & 0.7400 & 0.5536 & 0.726
\tabularnewline
\midrule
\begin{tabular}{@{}c@{}}
$\norm{\omega_{\mathcal{B}}(t_{f})-\text{\ensuremath{\omega_{\mathcal{B},f}}}}_2$ 
\\ ($\text{\textdegree}/U_T$) \end{tabular}
& 0.3748 & 4.9455 & 0.7976 & 4.072
\tabularnewline
\bottomrule
\end{tabular}
\end{table}

\subsection{Initialization performance }

Now we test the initialization performance of the proposed approach by the following manner: First, the policy trained by the proposed approach generates trajectories from the initial conditions in the test set. The PTR method then employs these trajectories as an initial guess for trajectory optimization. We choose the PTR method for trajectory optimization since the PTR method has been successfully used in the powered descent guidance applications \cite{szmuk2018successive,szmuk2017successive}. We compare the number of iterations and success rate of our method with that of the conventional initialization method that is the straight-line interpolation \cite{malyuta2021convex}. The termination parameters $\epsilon_{\nu}$ and $\epsilon_{tr}$ are $10^{-6}$ and $10^{-3}$, respectively, in (\ref{eq:PTR_termination}). Since the performance of the PTR method depends on the choice of parameters $w_{\nu}$ and $w_{tr}$ in (\ref{eq:ptrtr}), the experiments are carried out with 18 different combinations of these parameters for every initial states, then we record the minimum number of iterations. In the combinations of parameters, $w_{\nu}$ take a value in $\{10^3,10^4,10^5\}$ and $w_{tr}$ take a value among $\{10^{-2},10^{-1},1,10,10^1,10^2,10^3\}$. The maximum number of iterations in the PTR method is set as 50 since successful convergence of the PTR method usually requires no more than 50 iterations \cite{szmuk2017successive}. If the number of iterations exceeds this maximum, we consider that the optimization has failed to converge.

Table~\ref{tb:initialization} summarizes the initialization performance results. The proposed approach provides an initial guess with which the PTR method is able to converge for every initial state in the test set, while the straight-line interpolation has 7 failure cases. The table shows the mean, median, and standard deviation of the iteration count only for the successful cases. Our method shows much faster convergence performance than that of the straight-line interpolation.
\addtolength{\textheight}{0cm}

\begin{table}
\caption{Initialization Performance Comparison}
\label{tb:initialization}
\centering{}%
\begin{tabular}{ccc}
\toprule 
Property & Proposed & Straight-line \tabularnewline
\midrule 
Convergence success rate & 100\% (100/100) & 93\% (93/100)\tabularnewline
\midrule 
Mean & 2.20 & 8.55\tabularnewline
\midrule 
Median & 2 & 5\tabularnewline
\midrule 
Standard deviation & 0.40 & 7.90\tabularnewline
\bottomrule
\end{tabular}
\end{table}

\section{Conclusions}
\label{sec5}

    This paper has studied the problem of generating reliable initial guesses for trajectory optimization methods. 
    SCP-based guided policy search was applied to train the neural network that can generate trajectories with different initial conditions. To train the neural net policy, the proposed approach employed an interactive approach between convex optimization with penalized trust region and supervised learning. Numerical evaluations show that the neural network can generate good trajectory initial guesses from different initial conditions. Using these trajectories as the initial trajectory guesses is beneficial for the fast and reliable convergence of trajectory optimization methods.

\bibliographystyle{IEEEtran}
\bibliography{root}

\end{document}